\documentclass[12pt,a4paper]{article}
\usepackage[T1]{fontenc}
\usepackage{mathptmx}
\usepackage{courier}
\usepackage[dvips]{graphicx}
\usepackage{psfrag}
\usepackage{url}

\def\ds{\displaystyle}

\title{Scaling Invariance and the Iterative Transformation Method for a Class of Parabolic Moving Boundary Problems}
\author{Riccardo Fazio \\
Department of Mathematics and Computer Science, \\
University of Messina, Viale F. Stagno D'Alcontres, 31 \\
98166 Messina, Italy \\
E-mail: rfazio@unime.it \\
Home-page: http://mat521.unime.it/fazio}
\date{\today}
\begin{document}
\maketitle

\vspace{-1cm}
\begin{abstract}
In this paper we apply a scaling invariance analysis to reduce a class of parabolic moving boundary problems to free boundary problems governed by ordinary differential equations.
As well known free boundary problems are always non-linear and, consequently, their numerical solution is often obtained iteratively.
Among the numerical methods, developed for the numerical solution of this kind of problems, we focus on the iterative transformation method that has been defined within scaling invariance theory.
Then, as illustrative examples, we solve two problems of interest in the applications.
The obtained numerical results are found in good agreement with exact or approximate ones. 
\end{abstract}
%\bigskip
%\bigskip


%\bigskip
%\bigskip

%\newpage

\section{Introduction.}
The main contribution of this paper is the development of a complete scaling invariance analysis for a class of parabolic moving boundary problems and apply an extended scaling transformation to solve the reduced free boundary problems by an initial value solver.
In fact, the scaling invariance analysis allows us to reduce the original class of problems to a class of free boundary problems governed by ordinary differential equations.
As pointed out first by Landau \cite{Landau:HCM:1950}, free or moving boundary problems are always non-linear, and therefore they are often solved numerically.
In fact, a superposition principle for the solution of these problems cannot be valid because a variation of the auxiliary data, initial or free boundary conditions, produces a change in the free or moving boundary that in turn changes the domain of existence of the solution.

Most of the existing references on this subject, such as the one in the book by Crank \cite{Crank:FMB:1984}, or the one by Tarzia \cite{Tarzia:BMF:1984}, are essentially devoted to the solution of the well known Stefan problem which is governed by the linear heat equation. 
For the numerical solution of such a problem  several different approaches have been conceived over the years. 
Among the most famous ones there are the front-tracking, the front-fixing, and the domain-fixing methods (see \cite[pp. 217-281]{Crank:FMB:1984}), as well as other finite-difference or finite element approaches (see, for instance, Meek and Norbury \cite{Meek:NMB:1984}, Bonnerot and Jamet \cite{Bonnerot:SOF:1974} or Asaithambi \cite{Asaithambi:GMS:1992}), or moving grid, level set, or phase field methods (see the review by Javierre \cite{Javierre:CNM:2004}). 
Unfortunately, the proposed methods are introduced in the case of linear parabolic partial differential equations and they are not easily extended to non-linear parabolic cases belonging also to the class of problem (\ref{model}).
In this paper, we propose a method to overcome these difficulties, provided that the problem is invariant with respect to a scaling group. 
To this end, we use the similarity approach, described in full details by Dresner \cite{Dresner:1983:SSN}, within Lie's group invariance theory (see Bluman and Cole \cite{Bluman:1974:SMD}, Dresner \cite{Dresner:1999:ALT}, Barenblatt \cite{Barenblatt:1996:SSI}, or Bluman and Kumei \cite{Bluman:1989:SDE}). 
As far as the performance of different methods is concerned,
the introductory remark in a survey paper by Fox \cite{Fox:WAB:1975}
is pertinent: ``Problems of the same general nature
can differ enough in detail to make a good method for one problem
less satisfactory and even mediocre for another almost similar
problem''.
This point of view justifies the development of so many different
numerical methods.

The application of scaling invariance to applied mathematics and numerical analysis
has been a fruitful research field for more than a century, see Fazio \cite{Fazio:1999:NAS}.
As far as numerical applications are concerned, the first numerical transformation method (TM) is due to T{\"o}pfer \cite{Topfer:1912:BAB}.
He solved non-iteratively, using a transformation of variables, the Blasius problem of boundary layer theory. 
Several problems in boundary-layer theory lack this kind of invariance and cannot be solved by non-iterative (I)TMs, see \cite[Chapters 7-9]{Na:1979:CME}. 
To overcome this drawback it is possible to define an iterative extension of the T{\"o}pfer's algorithm  \cite{Fazio:1990:SNA,Fazio:1994:NTM,Fazio:1994:FSEb,Fazio:1996:NAN}.
Numerical solution of free boundary problems by means of TMs, developed within scaling invariance theory, are considered in \cite{Fazio:1990:SNA,Fazio:1991:ITM,Fazio:1998:SAN}.
Here, we explain in full details the definition of the ITM. 
Moreover, in order to show the validity of the proposed approach we solve two relevant problems:
the single phase Stefan's problem \cite{Stefan:UEP:1889,Crank:FMB:1984,Hill:ODS:1987}, and a problem describing the spreading of a viscous fluid above a smooth horizontal surface \cite{Meek:NMB:1984}.
For the former we compare the results obtained by the proposed ITM with those given by an asymptotic analysis.
As far as the second problem is concerned, the exact similarity solution is available, and therefore we are able to present a direct test for the obtained numerical results.

\section{Scaling invariance}
We consider the following class of moving 
boundary problems of the parabolic type
\begin{eqnarray}
& {\ds \frac{\partial u}{\partial t}} = 
{\ds \frac{\partial}{\partial x}\left[u^n \frac{\partial u}{\partial x} \right]} \quad \mbox{on} \quad t > 0, \ \ 0 < x < x_w(t) \ , \nonumber \\
& u(x, 0) = 0 \ , \qquad x_w(0) = 0 \ ,  \nonumber \\
& u(0, t) = A t^{\alpha} \ , \quad \left[\ \mbox{aut} \ {\ds \frac{\partial u}{\partial x}(0,t) = B t^{\beta}} \ \right] \ , \label{model} \\
& u(x_w(t), t) = p {\ds\left(t,x_w(t),\frac{dx_w}{dt}(t)\right)} \ , \nonumber  \\
& {\ds \frac{\partial u}{\partial x}}(x_w(t), t) = 
q {\ds\left(t,x_w(t),\frac{dx_w}{dt}(t),u(x_w(t), t)\right)} \ , \nonumber 
\end{eqnarray}
where $n$, $A$, $\alpha$, $B$, and $\beta$ are constants, $x$ and $t$ represent time and space respectively, $u(x,t)$ is the field variable, $p(\cdot, \cdot, \cdot)$, and $q(\cdot, \cdot, \cdot, \cdot)$ are given functions of their arguments, and 
$x_w(t)$ is the unknown moving boundary. 

As mentioned before, the problem (\ref{model}) is non-linear 
because $ x_w(t) $ depends on the initial and boundary data
so that a superposition principle cannot be valid (that was pointed out by Landau \cite{Landau:HCM:1950}).
As a consequence, obtaining analytical solutions for problems  
belonging to the class (\ref{model}) is a difficult task (see \cite[pp. 101-139]{Crank:FMB:1984}).

\begin{figure}[p]
\centering
%\psfragscanon 
\psfrag{xt}[][]{$ O $} 
\psfrag{x}[][]{$ x $} 
\psfrag{x1}[][]{$ x_1(t) = \eta_1 \; t^{1/\gamma} $} 
\psfrag{xw}[][]{$ x_w(t) = \eta_w \; t^{1/\gamma} $} 
\psfrag{e}[][]{$ \eta $} 
\psfrag{e1}[][]{$ \eta_1 $} 
\psfrag{ew}[][]{$ \eta_w $} 
\psfrag{t}[][]{$ t $} 
\includegraphics[width=.7\textwidth]{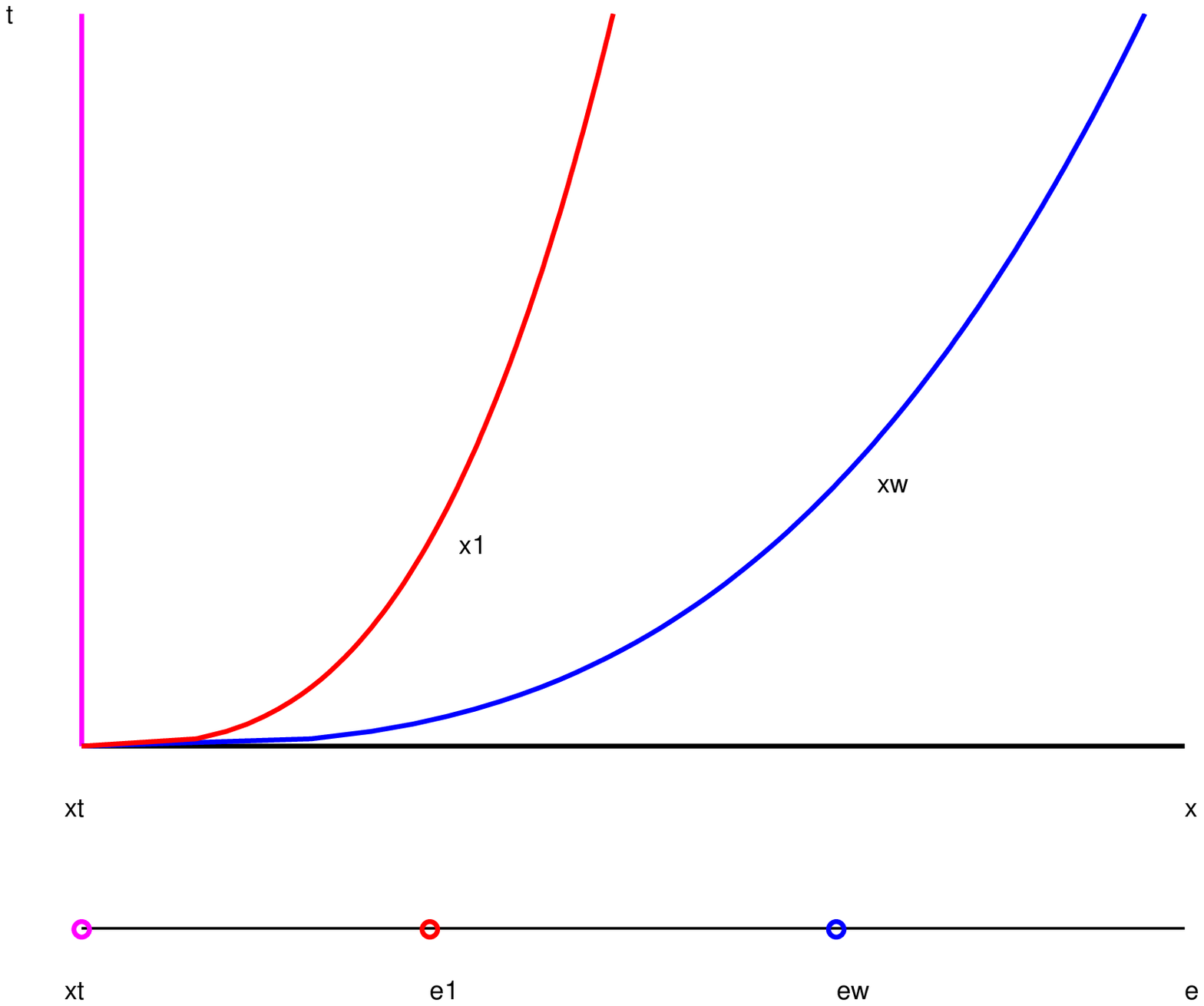}
\caption{The similarity map. Note that the similarity lines $ x = 0 $ and $ x_w (t) $ are mapped to $ \eta = 0 $ and $ \eta_w $, respectively.}
\label{fig:map}
\end{figure}

The following scaling group
\begin{equation}\label{sgroup}
x^* = \lambda x \ , \qquad {x_w}^* = \lambda x_w \ , \qquad 
t^* = \lambda^{\gamma} t \ , \qquad u^* = \lambda^{\alpha \gamma} u \ ,
\end{equation}
where  $\lambda $ is the (positive) group parameter, leaves the problem (\ref{model}) invariant provided that
\begin{eqnarray}
& \gamma = {\ds \frac{2}{n \alpha + 1}} \quad \left[\ \mbox{aut} \ {\ds \alpha = \frac{\beta + 1}{2 - n - n \beta}} \quad \mbox{with} \ {\ds \beta \ne \frac{2}{n} - 1} \ \right]\nonumber \\
& p(\cdot, \cdot, \cdot) = {\ds t^{\alpha} P\left(x_w(t) t^{-1/\gamma},\frac{dx_w}{dt}(t) t^{(\gamma - 1)/\gamma}\right)} \ , \label{eq:vincoli} \\
& q(\cdot, \cdot, \cdot, \cdot) = {\ds t^{(\alpha \gamma-1)/\gamma} Q\left(x_w(t) t^{-1/\gamma},\frac{dx_w}{dt}(t) t^{(\gamma-1)/\gamma},u(x_w(t), t)t^{-\alpha}\right)} \ . \nonumber
\end{eqnarray}
As a consequence, we can introduce the similarity variables as follows:
\begin{equation}\label{simvar}
\eta = x t^{-1/\gamma} \ , \qquad \eta_w = x_w(t) t^{-1/\gamma} \ , \qquad U(\eta) = t^{-\alpha} u(x,t) \ .
\end{equation}
By using (\ref{simvar}), we see that the model problem (\ref{model}) reduces to
\begin{eqnarray}
& {\ds \frac{d^2 U}{d\eta^2} + n U^{-1} \left(\frac{dU}{d\eta}\right)^2 + \frac{1}{\gamma} \eta U^{-n} \frac{dU}{d\eta}-\alpha U^{1-n}} = 0 \ , \nonumber \\
& U(0) = A \ , \quad \left[\ \mbox{aut} \ {\ds \frac{dU}{d\eta}(0)} = B \ \right] \ , \label{ODEs} \\
& U(\eta_w) = P\left(\eta_w, \eta_w\right) \ , \quad {\ds \frac{dU}{d\eta}(\eta_w) = Q\left(\eta_w, \eta_w, U(\eta_w)\right) } \ , \nonumber
\end{eqnarray}
where $\eta_w$ is the unknown free boundary for the ordinary differential problem (\ref{ODEs}).
We notice that, for any functional form of $P(\cdot, \cdot)$ and $Q(\cdot, \cdot, \cdot)$, the free boundary conditions in (\ref{ODEs}) depend only on $\eta_w$. 
Figure \ref{fig:map} shows the map given by the similarity variables (\ref{simvar}).

\section{The numerical method}
Here we consider the class of free boundary problems
\begin{eqnarray}
& {\ds \frac{d^2 w}{dz^2} = f\left(z,w,\frac{dw}{dz}\right)} \ , \nonumber \\[-1.5ex]
&  \label{Free} \\[-1.5ex]
& g\left(w(0), {\ds \frac{dw}{dz}}(0)\right) = C \ , \quad w(s) = j(s) \ , \quad {\ds \frac{w}{dz}}(s) = \ell(s) \ , \nonumber
\end{eqnarray}
where $w$ and $z$ are the field and independent variable, respectively, $f(\cdot,\cdot,\cdot)$, $g(\cdot,\cdot)$, $j(\cdot)$, and $\ell(\cdot)$ are given functions of their variables, $C$ is a given constant and $s$ represents the unknown free boundary. 
Of course, it is a simple matter to verify that (\ref{ODEs}) belongs to (\ref{Free}) for appropriate choices of $f$, $g$, $j$, and $\ell$.
A free boundary problem belonging to (\ref{Free}) can be solved by the ITM defined by the following steps.
\begin{itemize}
\item[-)] First we introduce an extended problem, namely:
\begin{eqnarray}
& {\ds \frac{d^2 w}{dz^2} = h^{(1-2\delta)/\sigma} f\left(h^{-\delta/\sigma} z,h^{-1/\sigma} w,h^{(\delta-1)/\sigma} \frac{dw}{dz}\right)} \ , \nonumber \\
& h^{1/\sigma} g\left(h^{-1/\sigma} w(0),h^{(\delta-1)/\sigma} {\ds \frac{dw}{dz}}(0)\right) = C  \ , \label{ExFree} \\
& w(s) = h^{1/\sigma} j(h^{-\delta/\sigma} s) \ , \quad {\ds \frac{w}{dz}}(h^{-\delta/\sigma} s) = h^{(1-\delta)/\sigma} \ell(h^{-\delta/\sigma} s) \ , \nonumber
\end{eqnarray}
where $C \ne 0$, and $h$ is a parameter.
A constructive characterization of (\ref{ExFree}), within similarity analysis, is given in \cite{Fazio:1997:NTE}.
Let us remark that the free boundary problem (\ref{Free}) is recovered from the extended problem (\ref{ExFree}) by setting $h=1$.
Moreover, the extended problem (\ref{ExFree}) is partially invariant with respect to the extended scaling group
\begin{equation}\label{Exgroup}
z^* = \omega^{\delta} z \ , \qquad s^* = \omega^{\delta} s \ , \qquad 
w^* = \omega w \ , \qquad h^* = \omega^{\sigma} h \ ,
\end{equation}
where $\omega$ is the (positive) group parameter, while $\delta$ and $\sigma$ are constants related to the particular problem under study.
We notice that the governing equations and the two free boundary conditions are invariant, but the condition at $z=0$ is not invariant.
\item[-)]
Given the values of $\delta$, $\sigma$ and $h^*$, we fix a value of $s^*$ greater than zero, and integrate (\ref{ExFree}), written in the starred variables, inwards in $[0, s^*]$ to compute approximate values of $w^*(0)$ and ${\ds \frac{dw^*}{dz^*}(0)}$ in order to get
\begin{equation}\label{omega}
\omega = {\ds \frac{h^{*1/\sigma} g\left(h^{*-1/\sigma} w^*(0),h^{*(\delta-1)/\sigma} {\ds \frac{dw^*}{dz^*}}(0)\right)}{C}} \ ,
\end{equation}
and, by using the scaling invariance, the corresponding value of $h = \omega^{-\sigma} h^*$.
\item[-)]
We get a solution of the original free boundary problem (\ref{Free}) when we find a value of $h^*$ that transforms to $h=1$.
This is equivalent to find a zero of the so called transformation function
\begin{equation}\label{Gamma}
\Gamma(h^*) = \omega^{-\sigma} h^* - 1 \ .
\end{equation}
To this end we can apply a root finder or a bracketing method.
The values of interest are defined by the scaling relations:
\begin{equation}\label{values}
s = \omega^{-\delta} s^* \ , \quad w(0) = \omega^{-1} w^*(0) \ , \quad 
{\ds \frac{dw}{dz}(0)  = \omega^{\delta - 1} \frac{dw^*}{dz^*}(0)} \ .
\end{equation}
Within the iteration we define the sequences $h^*_j$ and $s_j$ for $j = 0, 1, 2, \dots$.
If $\Gamma(h^*_j)$ tends to zero as $j$ goes to infinity, then $s_j$ goes to the correct free boundary value $s$ in the same limit.
\end{itemize}

For both examples reported in the next section we used, as initial value solver, the classical fourth order Runge-Kutta scheme and the secant method as root finder.
The convergence criterion for the secant method, was given by
\begin{equation}\label{eq:Tol}
|\Gamma(h_j^*) | \le \mbox{\rm Tol} \ , \quad \mbox{and} \quad |s_j - s_{j-1}| \le \mbox{\rm Tol} \ ,
\end{equation}
where $\mbox{\rm Tol} = 1\mbox{D}-06$. 
Here and in the following, the $\mbox{D}$ notation indicate a double precision arithmetic.

\section{Two applications}
As a first example, we consider the classical parabolic moving 
boundary problem: the celebrated Stefan's problem \cite[pp. 2-4, p. 9]{Crank:FMB:1984}.
This problem, in non-dimensional variables, can be written as follows
\begin{eqnarray}
& {\ds \frac{\partial u}{\partial t}} = 
{\ds \frac{\partial^2 u}{\partial x^2}} \quad \mbox{on} \quad t > 0 \ , \ \ 0 < x < x_w(t) \ , \nonumber \\
&  u(x, 0) = 0 \ , \qquad x_w(0) = 0 \ ,  \nonumber \\
& u(0, t) = 1 \ , \label{Stefan} \\
& u(x_w(t), t) = 0 \ , \nonumber \\
& {\ds \frac{\partial u}{\partial x}}(x_w(t), t) = 
- S {\ds \frac{dx_w}{dt}(t)} \ . \nonumber 
\end{eqnarray}
The last boundary condition, involving the free boundary velocity, is known as the Stefan's equation, and $1/S$ is the Stefan's number.
The scaling invariance analysis of section 2 is valid for (\ref{Stefan}) by setting $n=0$, $A=1$, $\alpha=0$, $\gamma=2$, $P(\cdot, \cdot)=0$, and $Q(\cdot, \cdot, \cdot)=-S t^{1/2} {\ds \frac{dx_w}{dt}}$.
Therefore, by introducing the similarity variables (\ref{simvar}), with $\alpha=0$ and $\gamma=2$, the problem (\ref{Stefan}) reduces to
\begin{eqnarray}
& {\ds \frac{d^2 U}{d\eta^2} + \frac{1}{2} \eta \frac{dU}{d\eta}} = 0 \ , \nonumber \\
& U(0) = 1 \ , \label{ODE1} \\
& U(\eta_w) = 0 \ , \quad {\ds \frac{dU}{d\eta}(\eta_w) = - \frac{1}{2} S \eta_w} \ . \nonumber
\end{eqnarray}
 Newmann's solution can be easily expressed by using the special error function 
$\mbox{\rm erf}(\cdot)$, see \cite{Abramowitz:HMF:1964},
\begin{equation}\label{Neumann}
U(\eta) = 1 - {\ds \frac{\mbox{\rm erf}(\eta/2)}{\mbox{\rm erf}(\eta_w/2)}} \ .
\end{equation}
The free boundary, $\eta_w$, is the positive real root of the equation
\begin{equation}\label{eq:root}
\pi^{1/2} S \eta_w \exp({\eta_w}^2/4) \mbox{\rm erf}(\eta_w/2) - 2 = 0 \ .
\end{equation}
This equation is transcendental and its solution allows us to obtain the exact moving boundary solution. 
The solution of (\ref{eq:root}) is unique ,and this implies the uniqueness of the similarity solution. 
As a consequence, the Stefan's problem admits only one solution.
Approximate values of $\eta_w$, for different values of $S$, can be found by a standard error function asymptotic expansion, see \cite{Hill:ODS:1987}.

\begin{figure}[p]
	\centering
\psfrag{e}[][]{$ \eta $} 
\psfrag{c1}[][]{$ U(\eta) $} 
\psfrag{c2}[][]{${\displaystyle \frac{dU}{d\eta}}(\eta) $} 

\includegraphics[width=.7\textwidth]{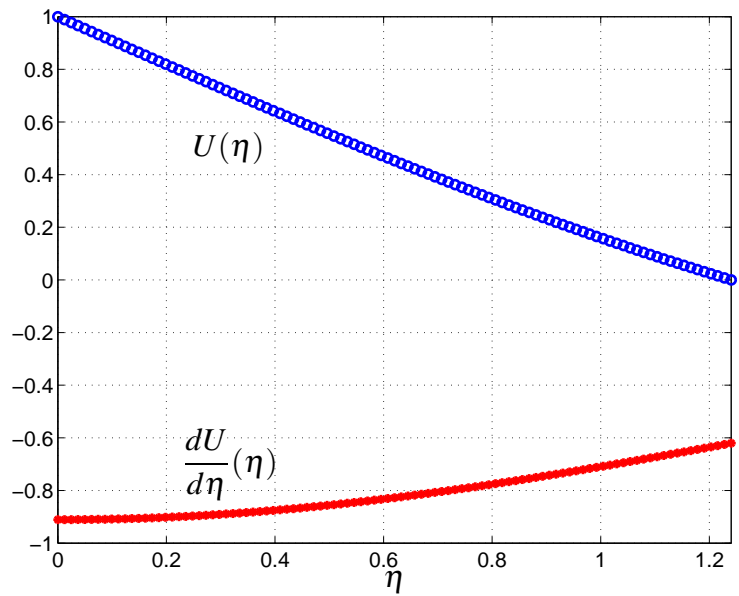}
\caption{Numerical solution for Stefan's problem with $S = 1$.
We used the results reported in table \ref{tab:Stefan} and the Runge-Kutta method with one hundred of steps.} 
	\label{fig:Stefan}
\end{figure}

In the following we want to find approximate values of $\eta_w$ by using the ITM.
First we have to introduce an extended problem
\begin{eqnarray}
& {\ds \frac{d^2 U}{d\eta^2} + \frac{h^{1/2}}{2} \eta \frac{dU}{d\eta}} = 0 \ , \nonumber \\
& U(0) = 1 \ , \label{ExODE1} \\
& U(\eta_w) = 0 \ , \quad {\ds \frac{dU}{d\eta}(\eta_w) = - \frac{h^{3/4}}{2} S \eta_w} \ , \nonumber
\end{eqnarray}
and the extended scaling group 
\begin{equation}\label{Exgroup1}
\eta^* = \omega^{-1} \eta \ , \qquad {\eta_w}^* = \omega^{-1} \eta_w \ , \qquad 
U^* = \omega U \ , \qquad h^* = \omega^{4} h \ .
\end{equation}
So that, by using the scaling invariance, it follows that
\begin{equation}\label{values1}
\omega = U^*(0)  \ , \quad h = \omega^{-4} h^* \ , \quad 
{\ds \frac{dU}{d\eta}(0)  = \omega^{- 2} \frac{dU^*}{d\eta^*}(0)} \ , \quad \eta_w = \omega {\eta_w}^* \ .
\end{equation}
Let us remark that, due to the $h^{1/2}$ term in (\ref{ExODE1}), we are allowed to consider positive values of $h^*$ only.
In table \ref{tab:Stefan} we list numerical results obtained for several values of $S$.
Comparing the data reported in the last two columns of table \ref{tab:Stefan}, we see that there is a good agreement between the results obtained by the present approach and those obtained by the asymptotic one. 
For the sake of brevity we omitted to list the intermediate data for the reported iterations.
The secant method always verified the convergence criterion (\ref{eq:Tol}) in few iterations.
Figure \ref{fig:Stefan} shows a sample numerical solution.

\begin{figure}[p]
	\centering
\psfrag{e}[][]{$ \eta $} 
\psfrag{c1}[][]{$ U(\eta) $} 
\psfrag{c2}[][]{${\displaystyle \frac{dU}{d\eta}}(\eta) $} 

\includegraphics[width=.7\textwidth]{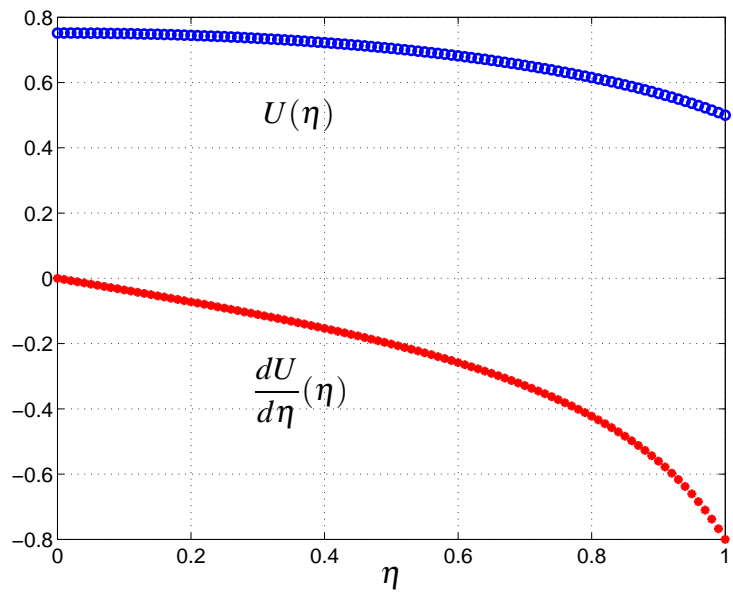}
\caption{Numerical solution for the problem (\ref{ODE}) with $L = 1/2$ and $H = -1/2$.
We notice that the exact similarity solution verifies the free boundary conditions: $U(\eta_w) = 1/2$ and 
$ {\ds \frac{dU}{d\eta}}(\eta_w) = -4/5$.
We used the results reported in the last line of table \ref{tab:Meek} and the Runge-Kutta method with one hundred of steps.} 
	\label{fig:Meek}
\end{figure}

As a second example, we consider a problem describing the spreading of a viscous fluid, such as a teacle, under the action of gravity above a smooth horizontal surface, see \cite{Meek:NMB:1984}.
Let us introduce the mathematical model
\begin{eqnarray}
& {\ds \frac{\partial u}{\partial t}} = 
{\ds \frac{\partial}{\partial x}\left[u^3 \frac{\partial u}{\partial x} \right]} \quad \mbox{on} \quad t > 0 \ , \ \ 0 < x < x_w(t) \ , \nonumber \\
&   u(x, 0) = 0 \ , \qquad x_w(0) = 0 \ ,  \nonumber \\
& {\ds \frac{\partial u}{\partial x}(0, t) = 0} \ , \label{Meek} \\
& u(x_w(t), t) = H t^{-1/5} \ , \nonumber \\
& {\ds \frac{\partial u}{\partial x}}(x_w(t), t) = 
L {x_w}^{-1} {\ds \frac{dx_w}{dt}(t)} u(x_w(t),t)^{-3} \ , \nonumber 
\end{eqnarray}
where $H$ and $L$ are given constants.
The scaling invariance analysis of section 2 is valid for (\ref{Meek}) by setting $n=3$, $B=0$, $\beta=-1/5$, $\gamma=5$, $P(\cdot, \cdot)=H$, and $Q(\cdot, \cdot, \cdot)= L t^{2/5} {x_w}^{-1} {\ds \frac{dx_w}{dt}(t)} u(x_w(t),t)^{-3}$.
Therefore, by using the similarity variables (\ref{simvar}), the free boundary problem for (\ref{Meek}) is given by
\begin{eqnarray}
& {\ds \frac{d^2 U}{d\eta^2} + 3 U^{-1} \left(\frac{dU}{d\eta}\right)^2 + \frac{1}{5} \eta U^{-3} \frac{dU}{d\eta}+\frac{1}{5} U^{-2}} = 0 \ , \nonumber \\
& {\ds \frac{dU}{d\eta}(0)} = 0 \ , \label{ODE} \\
& U(\eta_w) = H \ , \quad {\ds \frac{dU}{d\eta}(\eta_w) = \frac{L}{5H^3}} \ . \nonumber
\end{eqnarray}
We notice that the boundary condition at $\eta=0$ is homogeneous, that is $C=0$ in (\ref{Free}), and, as a consequence, in order to apply the ITM, we have to introduce a new dependent variable, namely
\begin{equation}\label{eq:newvar}
V(\eta) = U(\eta) + \eta \ .
\end{equation}
Hence, the free boundary (\ref{ODE}) becomes
\begin{eqnarray}
& {\ds \frac{d^2 V}{d\eta^2} + 3 (V-\eta)^{-1} \left(\frac{dV}{d\eta}-1\right)^2} + {\ds \frac{1}{5} \eta (V-\eta)^{-3} \left(\frac{dV}{d\eta}-1\right)+\frac{1}{5} (V-\eta)^{-2}} = 0 \ , \nonumber \\
& {\ds \frac{dV}{d\eta}(0)} = 1 \ , \label{ODE2} \\
& V(\eta_w) = H + \eta_w \ , \quad {\ds \frac{dV}{d\eta}(\eta_w) = \frac{L}{5H^3} + 1} \ . \nonumber
\end{eqnarray}
In order to apply the ITM, as a first step, we have to introduce an extended problem
\begin{eqnarray}
& {\ds \frac{d^2 V}{d\eta^2} + 3 (V-h^{1/2}\eta)^{-1} \left(\frac{dV}{d\eta}-h^{1/2}\right)^2} \nonumber \\
& \qquad + {\ds \frac{h^2}{5} \eta (V-h^{1/2}\eta)^{-3} \left(\frac{dV}{d\eta}-h^{1/2}\right)+\frac{h^2}{5} (V-h^{1/2}\eta)^{-2}} = 0 \ , \nonumber \\
& {\ds \frac{dV}{d\eta}(0)} = 1 \ , \label{ExODE2} \\
& V(\eta_w) = h H + h^{1/2} \eta_w \ , \quad {\ds \frac{dV}{d\eta}(\eta_w) = h^{1/2} \left(\frac{L}{5H^3} + 1\right)} \ . \nonumber
\end{eqnarray}
and the extended scaling group 
\begin{equation}\label{Exgroup2}
\eta^* = \omega^{1/2} \eta \ , \qquad {\eta_w}^* = \omega^{1/2} \eta_w \ , \qquad V^* = \omega V \ , \qquad h^* = \omega h \ . 
\end{equation}
Applying the scaling invariance properties we get
\begin{equation}\label{values2}
\omega = \left({\ds \frac{dU^*}{d\eta^*}(0)}\right)^2 \ , \quad h = \omega^{-1} h^* \ , \quad 
U(0) = \omega^{-1} V^*(0) \ ,  \quad {\eta_w} = \omega^{-1/2} {\eta_w}^* \ .
\end{equation}
Now, let us compare our numerical results with an exact similarity solution available in literature for specific values of the parameters.
In fact, in the particular case where $H= 1/2$ and $L=-1/2$ the exact solution
\begin{equation}\label{exact2}
U(\eta) = \left[{\ds \frac{3}{10} \left(\frac{5}{12}+{\eta_w}^2-\eta^2\right)}\right]^{1/3} \ , \qquad \eta_w = 1 \ ,
\end{equation}
was quoted in \cite{Meek:NMB:1984}.
In this case, from (\ref{exact2}) we get $U(0) \approx 0.751847$, that is correct to six decimal places.
Table \ref{tab:Meek} lists some numerical results obtained by the present approach for this particular choice of the parameter values.
At it is easily seen, the reported values for $U(0)$ and $\eta_w$ are in a very good agreement with the exact solution (\ref{exact2}).
Moreover, from the same table, we realize that comparable results must be found with different choices of ${\eta_w}^*$ and that, as far as the use of the secant method is concerned, we do not need to bracket the root of the transformation function. 
Figure \ref{fig:Meek} shows the obtained numerical solution.

As far as the model (\ref{Meek}) is concerned, the fluid flux at $x=0$ is given by
\begin{equation}\label{eq:flux}
u(0,t)^n {\ds \frac{\partial u}{\partial x}(0,t)} = B t^{\beta (n+1)} \left[U(0)\right]^n  {\ds \frac{dU}{d\eta}(0)}  \ .
\end{equation}
 However, the problem considered in \cite{Meek:NMB:1984} and studied herein for a viscous fluid has zero flux at $x=0$ because $B=0$.
On the other hand, the value of $U(0)$ is of interest because it defines the fluid height at $x=0$ according to
\begin{equation}\label{eq:height}
u(0,t) = t^{\beta} U(0)  \ ,
\end{equation}
where $\beta =-1/5$ for the considered fluid.

\section{Concluding remarks}
This paper makes two main contributions: the scaling invariance analysis for the class of parabolic moving boundary problems (\ref{model}) allows us to characterize those problems that can be reduced by similarity variables to free boundary problems governed by ordinary differential equations, and the definition of the ITM for the numerical solution of this second type of problems.
From a numerical viewpoint, it is simpler to solve a free boundary problem governed by an ordinary differential equation than to face a moving boundary problem governed by a parabolic partial differential equation.
The development of the ITM was motivated by the limitations of non-ITMs that sometimes can be applied due to the scaling invariance properties of the considered model, see \cite[Chapters 7-9]{Na:1979:CME}. 
The ITM has found application to the numerical solution of parabolic problems, either moving boundary problems or problems defined on infinite domains.  
In particular, in \cite{Fazio:2001:ITM} the ITM is used to solve the sequence of free boundary problems obtained by a semi-discretization of 1D parabolic moving boundary problems, and in \cite{Fazio:2010:MBF} a free boundary formulation for the reduced similarity models is used in order to propose a moving boundary formulation for parabolic problems on unbounded domains. 

 The present approach can be used for solving problems with different values of the parameters involving in
(\ref{model}) as well as different functional forms of 
\begin{eqnarray*}
&p {\ds\left(t,x_w(t),\frac{dx_w}{dt}(t)\right)} = {\ds t^{\alpha} P\left(x_w(t) t^{-1/\gamma},\frac{dx_w}{dt}(t) t^{(\gamma - 1)/\gamma}\right)} \\
&q {\ds\left(t,x_w(t),\frac{dx_w}{dt}(t),u(x_w(t), t)\right)} = \hspace{5cm}\\
& \qquad \qquad \qquad \qquad = {\ds t^{(\alpha \gamma-1)/\gamma} Q\left(x_w(t) t^{-1/\gamma},\frac{dx_w}{dt}(t) t^{(\gamma-1)/\gamma},u(x_w(t), t)t^{-\alpha}\right)} \ .
\end{eqnarray*}
Two problems, that have been defined in the applied sciences, were studied within the proposed framework.
The numerical results reported in the previous section clearly show the correctness and reliability of our approach.

\bigskip
\bigskip

\noindent
{\bf Acknowledgements.} This work was supported by the University of Messina.

%\begin{table}[!htb]
\begin{table}[p]
\renewcommand\arraystretch{1.1}
	\centering
		\begin{tabular}{rrr@{.}lccc}
\hline \\[-3ex]
\multicolumn{6}{c}%
%& 
%& \multicolumn{2}{c}%
%& \multicolumn{2}{c}%
%& \multicolumn{2}{c}%
{ITM}
&
{Asymptotic \cite{Hill:ODS:1987}}\\
{$S$} & {$j$}
& \multicolumn{2}{c}%
{${h_j^*}$}
& 
{${\displaystyle \frac{dU}{d\eta}(0)}$}
& 
{$\eta_w$} 
& 
{$\eta_w$}\\[2ex]
\hline
$0.1$ & $0$ & $600$ &  & & & \\
& $1$ & $700$ &  & & & \\
%& \vdots & \vdots &  & & & \\
& $10$ & $639$ & $263216$ & $-0.610425$ & $2.514145$ & $2.513961$ \\
$0.5$ & $0$ & $100$ &  & & & \\
& $1$ & $150$ &  & & & \\
& $9$ & $105$ & $180667$ & $-0.760017$ & $1.601231$ & $1.601187$ \\
$1.0$ & $0$ & $30$ &  & & & \\
& $1$ & $40$ &  & & & \\
& $8$ & $37$ & $843777$ & $-0.910875$ & $1.240134$ & $1.240161$ \\
$5.0$ & $0$ & $3$ &  & & & \\
& $1$ & $2$ &  & & & \\
& $7$ & $2$ & $256999$ & $-1.683000$ & $0.612848$ & $0.612864$ \\
$10.0$ & $0$ & $1$ &  & & & \\
& $1$ & $0$ & $5$  & & & \\
& $8$ & $0$ & $599873$ & $-2.309323$ & $0.440033$ & $0.440000$ \\
$50.0$ & $0$ & $1$ & $\mbox{D}-03$  & & & \\
& $1$ & $1$ & $\mbox{D}-02$ & & & \\
& $11$ & $2$ & $53\mbox{D}-02$ & $-5.03323$ & $0.199338$ & $0.199499$ \\
\hline			
		\end{tabular}
	\caption{Iterations obtained with ${\eta_w}^*=0.5$ and a step size $\Delta \eta^* = -1\mbox{D}-03$.}
	\label{tab:Stefan}
\end{table}

%\begin{table}[!htb]
\begin{table}[p]
\renewcommand\arraystretch{1.1}
	\centering
		\begin{tabular}{crr@{.}lr@{.}lcc}
\hline \\[-3ex]
{${\eta_w}^*$} & {$j$}
& \multicolumn{2}{c}%
{${h_j^*}$}
& \multicolumn{2}{c}%
{$\Gamma(h_j^*)$}
& {$U(0)$}
& {$\eta_w$} 
\\[2ex]
\hline
$0.5$ & $0$ & $0$ & $5$ & $0$ & $177999$ & & \\
& $1$ & $0$ & $1$ & $-0$ & $198207$ & & \\
& $7$ & $0$ & $250158$ & $-3$ & $14\mbox{D}-08$ & $0.751803$ & $0.996840$ \\
$1.0$ & $0$ & $0$ & $5$ & $-0$ & $152895$ & & \\
& $1$ & $0$ & $1$ & $-0$ & $349655$ & & \\
& $7$ & $1$ & $00032$ & $-1$ & $78\mbox{D}-08$ & $0.751825$ & $0.999842$ \\
\hline			
		\end{tabular}
	\caption{In the iterations we used a step size $\Delta \eta^* = -5\mbox{D}-04$.}
	\label{tab:Meek}
\end{table}

\end{document}